\documentclass[12pt]{article}
\usepackage{bbm}
\usepackage{latexsym}
\usepackage{mathrsfs}
\usepackage{amsmath}
\usepackage{graphicx}
\usepackage{amssymb}
\usepackage{color}
\usepackage{amsthm}
\usepackage{cite}
\usepackage{indentfirst}
\usepackage{anysize}\marginsize{45mm}{45mm}{40mm}{50mm}
\usepackage{caption}
\usepackage{enumerate}

\renewcommand\thesection{\arabic{section}.\kern -.5em}

\newtheoremstyle{lemma}{\topsep}{\topsep}%
     {}
     {}
     {\bfseries}
     {}
     {0.1em}
     {\thmname{#1}\thmnumber{ #2}\thmnote{ #3}}
\theoremstyle{lemma}  

\newtheorem{theorem}{Theorem}              
\newtheorem{lemma}[theorem]{Lemma}

\newtheorem{definition}{Definition}

\numberwithin{equation}{section}

\usepackage{latexsym, bm}
\begin{document}
\title{Hamiltonian cycles of balanced hypercube with more faulty edges\thanks{This research was partially supported by the National
Natural Science Foundation of China (No. 11801061).}}

\author{Ting Lan and Huazhong L\"{u}\thanks{Corresponding author.}\\
{\small School of Mathematical Sciences,}\\
{\small University of Electronic Science and Technology of China,}\\
{\small Chengdu, Sichuan, 610054, PR China}\\
{\small}lvhz@uestc.edu.cn}
\date{}

\maketitle

\begin{abstract}
The balanced hypercube $BH_{n}$, a variant of the hypercube, is a novel interconnection network for massive parallel systems. It is known that the balanced hypercube remains Hamiltonian after deleting at most $4n-5$ faulty edges if each vertex is incident with at least two edges in the resulting graph for all $n\geq2$. In this paper, we show that there exists a fault-free Hamiltonian cycle in $BH_{n}$ for $n\ge 2$ with $\left | F \right |\le 5n-7$ if the degree of every vertex in $BH_{n}-F$ is at least two and there exists no $f_{4}$-cycles in $BH_{n}-F$, which improves some known results.

\vskip 0.1 in


\end{abstract}

\section{Introduction}
In parallel and distributed systems, interconnection network plays an important role in their performance and specific functions.  The topological structure of an interconnection network is usually modeled by a simple undirected graph, where vertices represent processors and edges represent communication link between processors. The $n$-dimensional hypercube $Q_n$ is one of the most efficient interconnection network models, which has numerous good properties such as high symmetry, large connectivity and recursive structure. The balanced hypercube, $BH_{n}$, is a variant of the hypercube, which not only retains most of the excellent properties of $Q_n$, but also possesses some outstanding features that $Q_n$ does not have. For instance, the diameter of $BH_{n}$ is less than the hypercube with the same number of vertices for odd $n$. Moreover, for each vertex of $BH_{n}$, there is a unique vertex (backup vertex) with the same neighborhood, implying that if a vertex fails then tasks running on it can be easily transferred to its backup vertex without modifying routing. With such good properties, $BH_{n}$ has received much attention. Embedding properties, especially Hamiltonian path and cycle embedding of $BH_{n}$ have been studied in \cite{Chen,Cheng,Hao,Li,Xu,Yang,Qingguo,Huazhong}. Reliability indices under different theoretical model of $BH_{n}$ have been determined in \cite{Lu,YangD,YangM}. Some other topological issues of $BH_{n}$ have been investigated in \cite{lhlx,Zhou,Huang}.

\vskip 0.0 in

A graph $G$ is called {\em Hamiltonian} if there is a Hamiltonian cycle in $G$. A bipartite graph $G=(X\cup Y,E)$ is {\em Hamiltonian laceable} if there is a Hamiltonian path between any two vertices $x$ and $y$ in different partite sets. When designing and selecting the desired topology for an interconnected network, the Hamiltonian property is an important indicator for evaluating network performance \cite{Leighton}. Hamiltonian properties of the hypercube and its variants, as well as other famous networks have been widely studied \cite{Fan,Hsieh,Hsieh2,Park,Tsai}.

\begin{figure}[h]
\centering
\includegraphics[width=65mm]{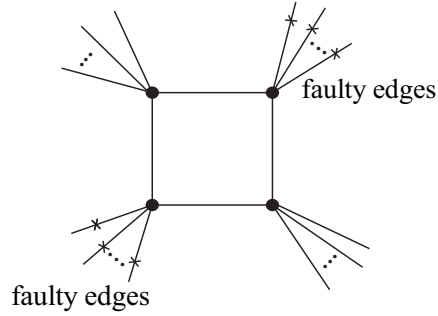}
\caption{$f_{4}$-$cycle$.} \label{f4cycle}
\end{figure}

In practical application, fault is inevitable in a large network, so it is necessary to consider fault-tolerance. A graph $G$ is {\em $k$ edge fault-tolerance Hamiltonian} if $G-F$ remains Hamiltonian for any $F\subset E(G)$ with $\left | F \right | \le k$. And a bipartite graph $G$ is {\em $k$ edge fault-tolerance Hamiltonian laceable} if $G-F$ remains Hamiltonian laceable for any $F\subset E(G)$ with $\left | F \right | \le k$. In \cite{Qingguo}, Zhou et al. proved that $BH_{n}$ is $2n-2$ edge fault-tolerant Hamiltonian laceable for $n\ge 2$. Recently, Li et al. \cite{Pingshan} further showed that there exists a Hamiltonian cycle in $BH_{n}-F$ for $n\ge 2$ with $\left | F \right |\le 4n-5$ if the degree of every vertex in $BH_{n}-F$ is at least two. They also showed that the upper bound $4n-5$ of the number of faulty edges can be tolerated is optimal. Actually, if a set $F$ of $4n-4$ edges is deleted from $BH_{n}$, then there exists no Hamiltonian cycles as a cycle of length four with a pair of nonadjacent vertices whose degrees are both two may occur in $BH_{n}-F$ (see Fig. \ref{f4cycle}). For convenience, we call such a cycle an {\em $f_{4}$-cycle}. Thus, a natural question arises: what is the upper bound of the number of faulty edges can be tolerated if $f_{4}$-cycle is forbidden in $BH_{n}-F$? Motivated by this, Liu and Wang \cite{Liu} showed that the hypercube still has a Hamiltonian cycle in $Q_{n}-F$ if $\delta(Q_{n}-F)\geq2$ ($|F|\leq3n-8$) and there exists no $f_{4}$-cycles.

\vskip 0.0 in

In this paper, we show that there exists a fault-free Hamiltonian cycle in $BH_{n}$ for $n\ge 2$ with $\left | F \right |\le 5n-7$ if the degree of every vertex in $BH_{n}-F$ is at least two and there exists no $f_{4}$-cycles in $BH_{n}-F$.

The rest of this paper is organized as follows. In Section 2, we give two equivalent definitions of balanced hypercubes and some useful properties of $BH_{n}$. In Section 3, we present the main result (Theorem \ref{5n-7}) of this paper. Finally, we conclude this paper in Section 4.

\section{Some definitions and lemmas}

Let $G=(V,E)$ be an undirected simple graph, where $V$ is the vertex set of $G$ and $E$ is the edge set of $G$. We write $e=uv$ for an edge $e$ with endpoints $u$ and $v$. Let $N_{G}(u)$ be the set of vertices adjacent to a vertex $u$ in $G$ and $e(u)$ be the set of edges incident with $u$. The degree of $u$, denoted by $d(u)$, is the number of edges incident with $u$. Let $\delta(G)$ be the minimum degree of $G$. A graph $G$ is bipartite if $V(G)$ is the union of two disjoint sets $X$, $Y$ such that every edge has one endpoint in $X$ and the other in $Y$.   A path $P$ from $v_{0}$ to $v_{n}$ is a sequence of distinct vertices $v_{0}v_{1}\dots v_{n}$ that consecutive vertices are adjacent, denoted by $\langle v_{0},P,v_{n} \rangle$. A cycle $C$ is a closed path $\langle v_{0},v_{1},\dots,v_{n},v_{0} \rangle$. A {\em Hamiltonian cycle} (resp. path) in a graph is a cycle (resp. path) that passes through every vertex of the graph.

\vskip 0.0 in

In the following, Wu and Huang \cite{Huang} gave the definitions of $BH_{n}$.

\begin{definition}{\bf .} An $n$-dimensional balanced hypercube $BH_{n}$
 consists of $2^{2n}$ vertices $(a_{0},$
$\ldots,a_{i-1},$ $a_{i},a_{i+1},\ldots,a_{n-1})$, where
$a_{i}\in\{0,1,2,3\}(0\leq i\leq n-1)$. An arbitrary vertex $(a_{0},\ldots,a_{i-1},$
$a_{i},a_{i+1},\ldots,a_{n-1})$ in $BH_{n}$ has the following $2n$ neighbors:

\begin{enumerate}
\item $((a_{0}+1)$ mod $4,a_{1},\ldots,a_{i-1},a_{i},a_{i+1},\ldots,a_{n-1})$,\\
      $((a_{0}-1)$ mod $4,a_{1},\ldots,a_{i-1},a_{i},a_{i+1},\ldots,a_{n-1})$, and
\item $((a_{0}+1)$ mod $4,a_{1},\ldots,a_{i-1},(a_{i}+(-1)^{a_{0}})$ mod $4,a_{i+1},\ldots,a_{n-1})$,\\
      $((a_{0}-1)$ mod $4,a_{1},\ldots,a_{i-1},(a_{i}+(-1)^{a_{0}})$ mod $
      4,a_{i+1},\ldots,a_{n-1})$.
\end{enumerate}
\end{definition}

\vskip 0.0 in

In $BH_{n}$, the first coordinate $a_{0}$ of $(a_{0},a_{1},\ldots,a_{n-1})$ is called {\em inner index}, and other coordinate $a_{i}$($1\le i\le n-1$) {\em i-dimensional index}. Let $uv$ be an edge of $BH_{n}$. If $u$ and $v$ differ only the inner index, then $uv$ is said to be a {\em 0-dimensional edge}. If $u$ and $v$  differ not only the inner index, but also $i$-dimensional index, then $uv$ is called an {\em i-dimensional edge}.

\vskip 0.0 in

\begin{definition}{\bf .} An $n$-dimensional balanced hypercube $BH_{n}$ can be recursively constructed as follows:

\begin{enumerate}

\item $BH_{1}$ is a cycle of length four whose vertices are labeled clockwise as $0$, $1$, $2$, $3$, respectively.
\item $BH_{n}$ is constructed from four disjoint $BH_{n-1}$s. These four $BH_{n-1}$s are labeled by $BH_{n-1}^{0}$, $BH_{n-1}^{1}$, $BH_{n-1}^{2}$, $BH_{n-1}^{3}$, where each vertex in $BH_{n-1}^{j}$($0\le j\le 3$) has $j$ attached as the new $n$-dimensional index. Each vertex $(a_{0},a_{1},\ldots,a_{n-2},j)$ has two extra neighbors:

(a). $((a_{0}+1)$ mod $4,a_{1},\ldots,a_{n-2},(j+1)$ mod $4)$ and
$((a_{0}-1)$ mod $4,a_{1},\ldots,a_{n-2},\\(j+1)$ mod $4)$, which are in $BH_{n-1}^{j+1}$ if $a_{0}$ is even;

(b). $((a_{0}+1)$ mod $4,a_{1},\ldots,a_{n-2},(j-1)$ mod $4)$ and
$((a_{0}-1)$ mod $4,a_{1},\ldots,a_{n-2},\\(j-1)$ mod $4)$, which are in $BH_{n-1}^{j-1}$ if $a_{0}$ is odd.

\end{enumerate}
\end{definition}

\vskip 0.0 in

\begin{figure}
\begin{minipage}[t]{0.5\linewidth}
\centering
\includegraphics[width=35mm]{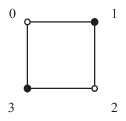}
\caption{$BH_{1}$.} \label{bh1}
\end{minipage}
\begin{minipage}[t]{0.5\linewidth}
\centering
\includegraphics[width=65mm]{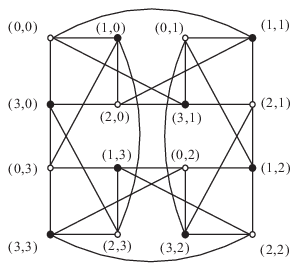}
\caption{$BH_{2}$.} \label{bh2}
\end{minipage}
\end{figure}

We illustrate $BH_{1}$ and $BH_{2}$ in Figs. \ref{bh1} and \ref{bh2}, respectively.

\vskip 0.0 in

Let $F$ be an edge subset and $D_{i}$ be the set of $i$-dimensional edges of $BH_{n}$, respectively, and let $F_{i}=F\cap D_{i}$ for $i=0,1,2,\dots,n-1$. Then, $BH_{n}$ can be divided into four subcubes $BH_{n-1}^{j,i}\cong BH_{n-1}$ by deleting $D_{i}$, where $j=0,1,2,3$ and $i\in\{0,1,2,\dots,n-1\}$. We use $BH_{n-1}^{j}$ to denote $BH_{n-1}^{j,n-1}$ for $j=0,1,2,3$ and let $F^{j}=F\cap E(BH^{j}_{n-1})$. For simplicity, we assume that $\left | F^{0} \right |=\max\left \{ \left | F^{j} \right | \mid  j=0,1,2,3  \right \}$.

\vskip 0.0 in

In the following, we make an convention that the vertices with superscript $j$ belong to $BH_{n-1}^{j,i}$, say $a^{j}$, $b^{j}$, $c^{j}$, $\dots \in V(BH^{j,i}_{n-1})$. An {\em $r$-edge} is an edge in $BH_{n-1}^{j,i}$ whose two endpoints are both incident with at least one $i$-dimensional nonfaulty edge, and if an edge in $BH_{n-1}^{j,i}$ has at least one endpoint incident with two $i$-dimensional faulty edges, we call it a {\em non-$r$-edge}. A {\em pivot vertex} is a vertex incident with exactly one non-faulty edge in $BH_{n-1}^{j,i}$.

\vskip 0.0 in

In what follows, we shall present some basic properties of balanced hypercube, which will be used later.

\begin{lemma}{\cite{Huang+Wu}\bf .}
$BH_{n}$ is bipartite.
\end{lemma}

\begin{lemma}{\cite{Zhou}\bf .}
$BH_{n}$ is vertex-transitive and edge-transitive.
\end{lemma}

\begin{theorem}{\label{2n-3}\cite{Huazhong}\bf .}
Let $F$ be a set of edges with $\left | F \right |\le 2n-3$ and let $\left \{ s_{1},s_{2} \right \}$ and $\left \{ t_{1},t_{2} \right \}$ be two sets of vertices in different partite sets of $BH_{n}$ for $n\ge 2$. Then $BH_{n}-F$ contains vertex-disjoint $s_{1},t_{1}$-path and $s_{2},t_{2}$-path that cover all vertices of it.
\end{theorem}

\begin{theorem}{\label{2n-2}\cite{Qingguo}\bf .}
$BH_{n}$ is $2n-2$ edge fault-tolerant Hamiltonian laceable for $n\ge 2$.
\end{theorem}

\begin{theorem}{\label{4n-5}\cite{Pingshan}\bf .}
Let $F$ be an arbitrary edge set with $\left | F \right |\le 4n-5$ and $\delta (BH_{n}-F)\ge 2$. Then each edge in $BH_{n}-F$ lies on a fault-free Hamiltonian cycle for all $n\ge 2$.
\end{theorem}

\begin{theorem}{\label{hyper}\cite{Zhang}\bf .}
$BH_{n}$ is hyper-Hamiltonian laceable for $n\ge 1$.
\end{theorem}

\section{Main Results}

\begin{lemma}{\label{8} \bf .}
Let $F$ be an edge subset of $BH_{3}$ with $\lvert F \rvert\le 8$. $BH_{3}-F$ is Hamiltonian if there exists no $f_{4}$-cycles in $BH_{3}-F$ and $\delta (BH_{3}-F)\ge 2$.

\vskip 0.05 in

\noindent {\bf Proof.} {\rm Clearly, we need only to consider that $\lvert F \rvert= 8$. Since $\sum_{i=0}^{2}\lvert F_{i} \rvert = \lvert F \rvert=8$, there is an integer $i\in \left \{0, 1, 2 \right \}$ such that $\left | F_{i} \right |\ge 3$. Without loss of generality, we assume $i=2$, that is $\left | F_{2} \right |\ge 3$. In the following, we decompose $BH_{3}$ into four subcubes $BH_{2}^{j}$ by deleting $D_{2}$, where $j=0, 1, 2, 3$. $\left | F^{j} \right | \le \sum_{j=0}^{3}\left | F^{j} \right |=\left | F \right |-\left | F_{2} \right |  \le 5$. We consider the following three cases in terms of the minimum degree of $BH_{2}^{j}-F^{j}$ and the cardinality of $F^{j}$, where $j=0,1,2,3$.

\vskip 0.02 in

\noindent {\bf Case 1.} $\delta (BH^{j}_{2}-F^{j})\ge 2$ for all $j\in \left \{0, 1, 2, 3\right \}$.

\vskip 0.02 in

\noindent {\bf Subcase 1.1.} $\left | F^{0} \right | \le3$.

\vskip 0.02 in

Since $\sum_{j=0}^{3}\left | F^{j} \right |=\left | F \right |-\left | F_{2} \right |   \le 5$ and $\left | F^{0} \right | \le3$, $\left | F^{j} \right | \le2$ for $j=1, 2, 3$. Obviously, there exists an $r$-edge $e_{0}=a^{0}b^{0}$ in $BH_{2}^{0}$. Let $a^{0} b^{1}$, $a^{3}b^{0}$ be the $2$-dimensional nonfaulty edges. By Theorem \ref{4n-5}, there is a fault-free Hamiltonian cycle $C_{0}$ in $BH_{2}^{0}-F^{0}$ that contains $e_{0}$. Let $H_{0}=C_{0}-e_{0}$ and let $a^{1}b^{2}$, $a^{2}b^{3}\in D_{2}-F_{2}$, by Theorem \ref{2n-2}, there is a Hamiltonian path $H_{j}$ in $BH^{j}_{2}-F^{j}$ that joins $a^{j}$ and $b^{j}$ for $j=1, 2, 3$. Thus, $\left \langle a^{0},b^{1},H_{1},a^{1},b^{2},H_{2},a^{2}, b^{3},H_{3},a^{3},b^{0}, H_{0}, a^{0}\right \rangle $ (see Fig. \ref{case11}) is a fault-free Hamiltonian cycle in $BH_{3}-F$.

\begin{figure}
\centering
\includegraphics[width=65mm]{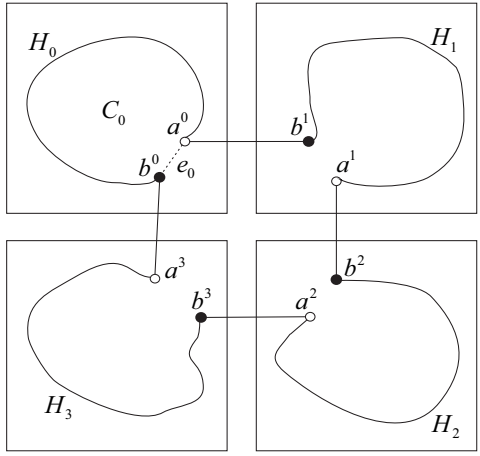}
\caption{Illustration for Subcase 1.1.} \label{case11}
\end{figure}

\vskip 0.02 in

\noindent {\bf Subcase 1.2.} $\left | F^{0} \right |=4$.

\vskip 0.02 in

Clearly, $3\le \left | F_{2} \right | \le 4$, and $\left | F^{j} \right | \le 1$ for each $j=1, 2, 3$. Since $\left | F_{2} \right | \le 4$, there is an edge $f_{0}=a^{0}b^{0}\in F^{0}$ adjacent to at least one $2$-dimensional nonfaulty edge, say $a^{0}b^{1}$. By Theorem \ref{4n-5}, there is a Hamiltonian cycle $C_{0}$ in $BH^{0}_{2}-F^{0}+f_{0}$ that contains $f_{0}$. Let $H_{0}=C_{0}-f_{0}$.

\vskip 0.02 in

If $f_{0}$ is an $r$-edge, assume that $a^{3}b^{0}\in D_{2}-F_{2}$. Let $a^{1}b^{2}$, $a^{2}b^{3}\in D_{2}-F_{2}$, by Theorem \ref{2n-2}, there is a Hamiltonian path $H_{j}$ in $BH^{j}_{2}-F^{j}$ that joins $a^{j}$ and $b^{j}$ for $j=1, 2, 3$. Thus, $\left \langle a^{0},b^{1},H_{1},a^{1},b^{2},H_{2},a^{2}, b^{3},H_{3},a^{3},b^{0}, H_{0}, a^{0}\right \rangle $ is a fault-free Hamiltonian cycle in $BH_{3}-F$.

\vskip 0.02 in

If $f_{0}$ is a non-$r$-edge, assume that $b^{0}$ is incident with two $2$-dimensional faulty edges. Let $c^{0}b^{0}\notin E(C_{0})$ with $c^{0}\in V(C_{0})$ and let $d^{0}$ and $e^{0}$ be the neighbors of $c^{0}$ satisfying $c^{0}d^{0}$, $c^{0}e^{0}\in E(C_{0})$. Let $H_{00}$ be the path from $a^{0}$ to $e^{0}$ on $C_{0}$ and let $H_{01}$ be the path from $d^{0}$ to $b^{0}$ on $C_{0}$, then we denote $C_{0}$ by $\langle a^{0},H_{00},e^{0},c^{0},d^{0},H_{01},b^{0},a^{0}\rangle $. Since $b^{0}$ is incident with two $2$-dimensional faulty edges and $\left | F_{2} \right | \le 4$, $d^{0}$ or $e^{0}$ is incident with at least one $2$-dimensional nonfaulty edge.

\vskip 0.02 in

(1) If $d^{0}$ is incident with a $2$-dimensional nonfaulty edge $a^{3}d^{0}$, let $a^{1}b^{2}$, $a^{2}b^{3}\in D_{2}-F_{2}$. Since $\left | F^{j} \right | \le 1$, by Theorem \ref{2n-2}, there is a Hamiltonian path $H_{j}$ in $BH^{j}_{2}-F^{j}$ that join $a^{j}$ and $b^{j}$ for $j=1,2,3$. Then $\langle a^{0},b^{1},H_{1},a^{1},b^{2},H_{2},a^{2},b^{3},H_{3},a^{3},d^{0},H_{01},b^{0},c^{0},e^{0},\\H_{00},a^{0}\rangle $ (see Fig. \ref{case12b}) is a fault-free Hamiltonian cycle in $BH_{3}-F$.

\vskip 0.02 in

(2) If $d^{0}$ is incident with two $2$-dimensional faulty edges, then there exists an edge $c^{3}e^{0}\in D_{2}-F_{2}$. We may assume that all edges in $F^{0}$ are non-$r$-edges, otherwise, we can replace $f_{0}$ by an $r$-edge in $F^{0}$. So all edges in $F^{0}$ are incident with $d^{0}$ or $b^{0}$. Since there exists no $f_{4}$-cycles in $BH_{3}-F$, the length of the cycle $C_{00}=\left \langle b^{0},c^{0},d^{0},H_{01},b^{0}\right \rangle $ is at least 6. Hence, there exists
an $r$-edge $g^{0}f^{0}$ on $C_{00}$. Let $H_{010}$ be the path from $d^{0}$ to $f^{0}$ on $C_{00}$ and let $H_{011}$ be the path from $g^{0}$ to $b^{0}$ on $C_{00}$. Then, we denote $C_{00}$ by $\left \langle b^{0},c^{0},d^{0},H_{010},f^{0},g^{0},H_{011},b^{0}\right \rangle $. Let $g^{0}d^{1}$, $a^{3}f^{0}\in D_{2}-F_{2}$ satisfying $c^{3}\ne a^{3}$ and $b^{1}\ne d^{1}$ and let $a^{j}b^{j+1}$, $c^{j}d^{j+1}\in D_{2}-F_{2}$ for $j=1,2$. By Theorem \ref{2n-3}, there are two vertex-disjoint paths $H_{j0}$ and $H_{j1}$ in $BH_{2}^{j}-F^{j}$ such that $H_{j0}$ joins $a^{j}$ and $b^{j}$ and $H_{j1}$ joins $c^{j}$ and $d^{j}$ for $i=1,2,3$. Then $\langle a^{0},b^{1},H_{10},a^{1},b^{2},H_{20},a^{2},b^{3},H_{30},a^{3},f^{0},H_{010},d^{0},c^{0},b^{0},H_{011},
g^{0},d^{1},H_{11},c^{1},d^{2},H_{21},c^{2},\\d^{3},H_{31},c^{3},e^{0},H_{00},a^{0}\rangle $ (see Fig. \ref{case12c}) is a fault-free Hamiltonian cycle in $BH_{3}-F$.

\begin{figure}
\begin{minipage}[t]{0.5\linewidth}
\centering
\includegraphics[width=65mm]{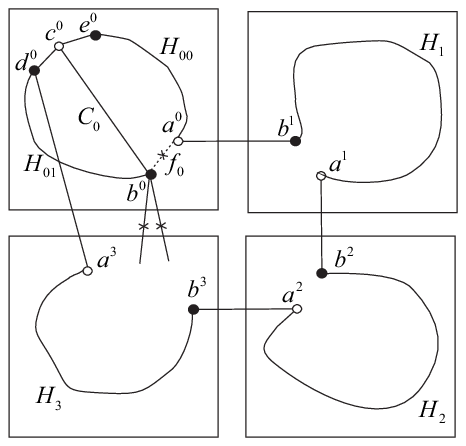}
\caption{Illustration for Subcase 1.2.} \label{case12b}
\end{minipage}
\begin{minipage}[t]{0.5\linewidth}
\centering
\includegraphics[width=65mm]{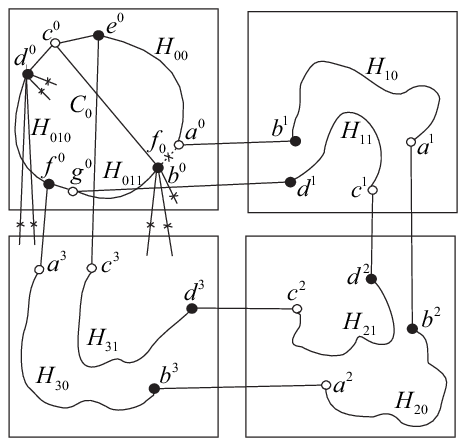}
\caption{Illustration for Subcase 1.2.} \label{case12c}
\end{minipage}
\end{figure}

\vskip 0.02 in

\noindent {\bf Subcase 1.3.} $\left | F^{0} \right |=5$.

\vskip 0.02 in

Then $\left | F^{1} \right |=\left | F^{2} \right |=\left | F^{3} \right |=0$ and $\left | F_{2} \right |=3$. Since $\left | F_{2} \right |=3$ and $\delta (BH^{0}_{2}-F^{0})\ge 2$, there exist two $r$-edges $f_{1}=a^{0}b^{0}$ and $f_{2}=c^{0}d^{0}$ in $F^{0}$. By Theorem \ref{4n-5}, there is a Hamiltonian cycle $C_{0}$ in $BH^{0}_{2}-F^{0}+f_{1}+f_{2}$ that contains $f_{1}$. If $f_{2}\notin E(C_{0})$, the situation is similar to the case that $f_{0}$ is $r$-edge in Subcase 1.2. So, we assume that $f_{2}\in E(C_{0})$.

\vskip 0.02 in

\noindent {\bf Subcase 1.3.1.} $f_{1}$ is adjacent to $f_{2}$ on $C_{0}$.

\vskip 0.02 in

Without loss of generality, we assume that $a^{0}=c^{0}$. Assume that $a^{3}b^{0}$, $c^{3}d^{0}$, $a^{0}b^{1}\in D_{2}-F_{2}$(possibly $a^{3}=c^{3}$). Let $a^{0}p^{0}$, $a^{0}q^{0}\notin E(C_{0})$ with $p^{0}$, $q^{0}\in V(C_{0})$ and let $x^{0}$, $y^{0}$ be the neighbors of $p^{0}$ satisfying $x^{0}p^{0}$, $y^{0}p^{0}\in E(C_{0})$ and $s^{0}$, $t^{0}$ be the neighbors of $q^{0}$ satisfying $s^{0}q^{0}$, $t^{0}q^{0}\in E(C_{0})$(possibly $y^{0}=s^{0}$). Since $\left | F_{2} \right |=3$, $s^{0}$, $t^{0}$, $x^{0}$ or $y^{0}$ is incident with a $2$-dimensional nonfaulty edge that is not incident with $b^{1}$. We may assume that $t^{0}d^{1}\in D_{2}-F_{2}$ with $d^{1}\ne b^{1}$. Let $H_{00}$ be the path from $b^{0}$ to $t^{0}$ on $C_{0}$ and $H_{01}$ be the path from $q^{0}$ to $d^{0}$ on $C_{0}$. We denote $C_{0}$ by $\left \langle a^{0},b^{0},H_{00},t^{0},q^{0},H_{01},d^{0},a^{0}\right \rangle$. Let $a^{j}b^{j+1}$, $c^{j}d^{j+1}\in D_{2}-F_{2}$ for $j=1,2$.

\vskip 0.02 in

If $c^{3}\ne a^{3}$, by Theorem \ref{2n-3}, there are two vertex-disjoint paths $H_{j0}$ and $H_{j1}$ in $BH_{2}^{j}-F^{j}$ such that $H_{j0}$ joins $a^{j}$ and $b^{j}$ and $H_{j1}$ joins $c^{j}$ and $d^{j}$ for $i=1,2,3$. Then $\langle a^{0}, b^{1}, H_{10}, a^{1}, b^{2}, H_{20}, a^{2}, b^{3}, H_{30}, a^{3}, b^{0}, H_{00}, t^{0}, d^{1}, H_{11}, c^{1}, d^{2},H_{21}, c^{2}, d^{3}, H_{31}, c^{3}, d^{0},\\ H_{01},q^{0}, a^{0}\rangle $ (see Fig. \ref{case131a}) is a fault-free Hamiltonian cycle in $BH_{3}-F$.

\vskip 0.02 in

If $c^{3}=a^{3}$, by Theorem \ref{hyper}, there is a Hamiltonian path $H_{3}$ in $BH^{3}_{2}-\left \{ a^{3} \right \}$ that joins $d^{3}$ and $b^{3}$.  And by Theorem \ref{2n-3}, there are two vertex-disjoint paths $H_{j0}$ and $H_{j1}$ in $BH_{2}^{j}-F^{j}$ such that $H_{j0}$ joins $a^{j}$ and $b^{j}$ and $H_{j1}$ joins $c^{j}$ and $d^{j}$ for $i=1,2$. Then, $ \langle a^{0},b^{1},H_{10},a^{1},b^{2},H_{20},a^{2},b^{3},H_{3},d^{3},c^{2}, H_{21}, d^{2},c^{1},H_{11},d^{1},t^{0},H_{00},b^{0},a^{3},d^{0},H_{01},q^{0},\\a^{0} \rangle $ (see Fig. \ref{case131b}) is a fault-free Hamiltonian cycle in $BH_{3}-F$.

\begin{figure}
\begin{minipage}[t]{0.5\linewidth}
\centering
\includegraphics[width=65mm]{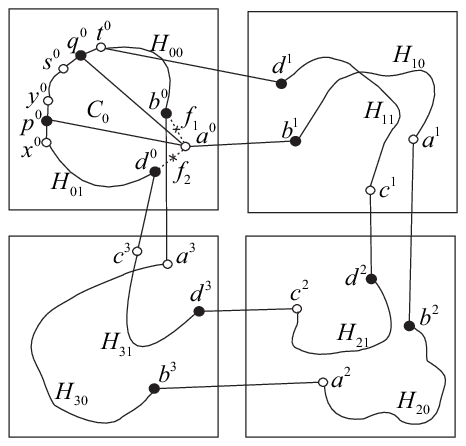}
\caption{Illustration for Subcase 1.3.1.} \label{case131a}
\end{minipage}
\begin{minipage}[t]{0.5\linewidth}
\centering
\includegraphics[width=65mm]{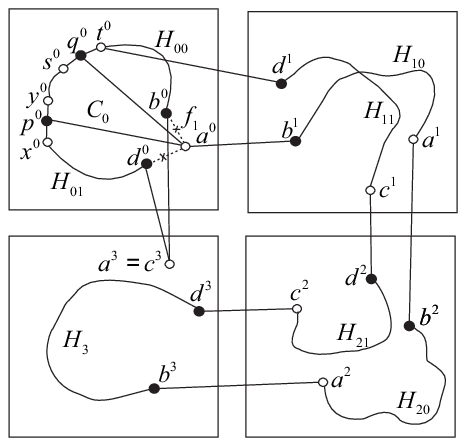}
\caption{Illustration for Subcase 1.3.1.} \label{case131b}
\end{minipage}
\end{figure}

\vskip 0.02 in

\noindent {\bf Subcase 1.3.2.} $f_{1}$ is nonadjacent to $f_{2}$ on $C_{0}$.

\vskip 0.02 in

We denote $C_{0}$ by $\left \langle a^{0},b^{0},H_{00},c^{0},d^{0},H_{01},a^{0}\right \rangle $. Since $\left | F_{2} \right |=3$ and $a^{0}b^{0}$, $c^{0}d^{0}$ are $r$-edges, $a^{0}$, $b^{0}$, $c^{0}$ and $d^{0}$ are incident with at least five $2$-dimensional nonfaulty edges in total, and each is incident with at least one $2$-dimensional nonfaulty edge. Then there must exist $c^{3}d^{0}$, $a^{3}b^{0}$, $a^{0}b^{1}$, $c^{0}d^{1}\in D_{2}-F_{2}$ satisfying $c^{3}\ne a^{3}$ or $b^{1}\ne d^{1}$. Without loss of generality, we assume that $b^{1}\ne d^{1}$. Let $a^{j}b^{j+1}$, $c^{j}d^{j+1}\in D_{2}-F_{2}$ for $j=1,2$.

\vskip 0.02 in

If $c^{3}\ne a^{3}$, by Theorem \ref{2n-3}, there are two vertex-disjoint paths $H_{j0}$ and $H_{j1}$ in $BH_{2}^{j}-F^{j}$ such that $H_{j0}$ joins $a^{j}$ and $b^{j}$ and $H_{j1}$ joins $c^{j}$ and $d^{j}$ for $i=1,2,3$. Then $\langle a^{0}, b^{1}, H_{10}, a^{1}, b^{2}, H_{20}, a^{2}, b^{3}, H_{30}, a^{3},b^{0}, H_{00}, c^{0}, d^{1}, H_{11}, c^{1}, d^{2},H_{21}, c^{2}, d^{3}, H_{31}, c^{3}, d^{0},\\ H_{01}, a^{0}\rangle $ (see Fig. \ref{case132a}) is a fault-free Hamiltonian cycle in $BH_{3}-F$.

\vskip 0.02 in

If $c^{3}=a^{3}$, by Theorem \ref{hyper}, there is a Hamiltonian path $H_{3}$ in $BH^{3}_{2}-\left \{ a^{3} \right \}$ that joins $d^{3}$ and $b^{3}$. And by Theorem \ref{2n-3}, there are two vertex-disjoint paths $H_{j0}$ and $H_{j1}$ in $BH_{2}^{j}-F^{j}$ such that $H_{j0}$ joins $a^{j}$ and $b^{j}$ and $H_{j1}$ joins $c^{j}$ and $d^{j}$ for $i=1,2$. Then, $\langle a^{0},b^{1},H_{10},a^{1},b^{2},H_{20},a^{2},b^{3},H_{3},d^{3},c^{2}, H_{21}, d^{2},c^{1},H_{11},d^{1},c^{0},H_{00},b^{0},a^{3},d^{0},H_{01},a^{0}\rangle$ \\ (see Fig. \ref{case132b}) is a fault-free Hamiltonian cycle in $BH_{3}-F$.

\begin{figure}
\begin{minipage}[t]{0.5\linewidth}
\centering
\includegraphics[width=65mm]{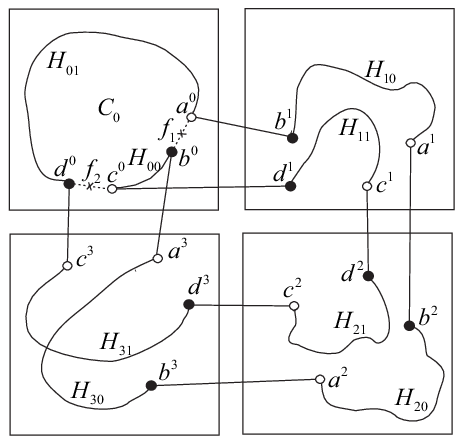}
\caption{Illustration for Subcase 1.3.2.} \label{case132a}
\end{minipage}
\begin{minipage}[t]{0.5\linewidth}
\centering
\includegraphics[width=65mm]{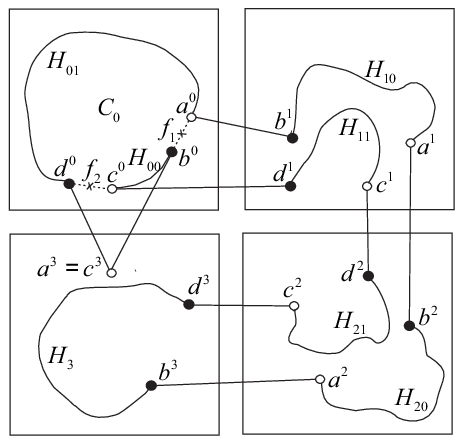}
\caption{Illustration for Subcase 1.3.2.} \label{case132b}
\end{minipage}
\end{figure}

\vskip 0.02 in

\noindent {\bf Case 2.} $\delta (BH^{j}_{2}-F^{j})=1$ for some $j\in \left \{0, 1, 2, 3\right \}$.

\vskip 0.02 in

\noindent {\bf Subcase 2.1.} $\left | F^{0} \right |=5$.

\vskip 0.02 in

Then $\left | F^{1} \right |=\left | F^{2} \right |=\left | F^{3} \right |=0$ and $\left | F_{2} \right |=3$. Since $\left | F^{0} \right |=5$, there are at most two pivot vertices in $BH_{2}^{0}$.

\vskip 0.02 in

\noindent {\bf Subcase 2.1.1.} There is exactly one pivot vertex $a^{0}$ in $BH_{2}^{0}$.

\begin{figure}
\begin{minipage}[t]{0.5\linewidth}
\centering
\includegraphics[width=65mm]{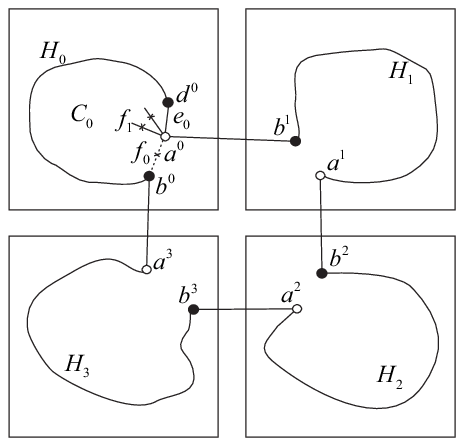}
\caption{Illustration for Subcase 2.1.1.} \label{case21}
\end{minipage}
\begin{minipage}[t]{0.5\linewidth}
\centering
\includegraphics[width=65mm]{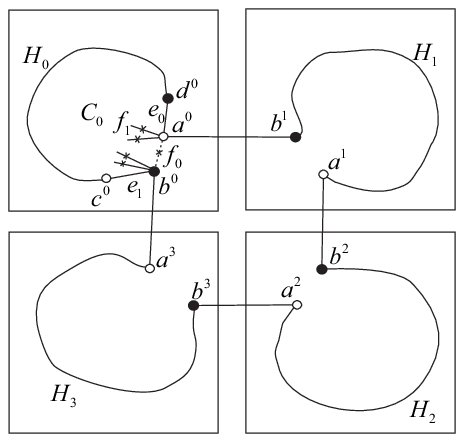}
\caption{Illustration for Subcase 2.1.2.} \label{case22}
\end{minipage}
\end{figure}

Let $f_{0}$ and $f_{1}$ be two $r$-edges in $e(a^{0})\cap F^{0}$ and let $e_{0}=a^{0}d^{0}$ be the nonfaulty edge in $e(a^{0})\cap E(BH_{2}^{0})$. By Theorem \ref{4n-5}, there is a Hamiltonian cycle $C_{0}$ in $BH^{0}_{2}-F^{0}+f_{0}+f_{1}$ that contains $e_{0}$. Notice also that one of $f_{0}$ and $f_{1}$ is on $C_{0}$, say $f_{0}=a^{0}b^{0}\in E(C_{0})$. Let $H_{0}=C_{0}-f_{0}$. Then we can choose $b^{0}a^{3}\in D_{2}-F_{2}$ and $a^{j}b^{j+1}\in D_{2}-F_{2}$ for $j=0,1,2$. By Theorem \ref{2n-2}, there is a Hamiltonian path $H_{j}$ in $BH^{j}_{2}-F^{j}$ that join $a^{j}$ and $b^{j}$ for $j=1,2,3$. Then $\left \langle a^{0},b^{1},H_{1},a^{1},b^{2},H_{2},a^{2},b^{3},H_{3},a^{3},b^{0},H_{0},a^{0}\right \rangle $ (see Fig. \ref{case21}) is a fault-free Hamiltonian cycle in $BH_{3}-F$.

\vskip 0.02 in

\noindent {\bf Subcase 2.1.2.} There are exactly two pivot vertices $a^{0}$ and $b^{0}$ in $BH_{2}^{0}$.

\vskip 0.02 in

Then $a^{0}$ and $b^{0}$ are endpoints of a faulty edge $f_{0}=a^{0}b^{0}\in F^{0}$. Let $f_{1}\in e(a^{0})\cap F^{0}-f_{0}$ and $e_{0}=a^{0}d^{0}\in E(BH^{0}_{2})-F^{0}$. By Theorem \ref{4n-5}, there exists a Hamiltonian cycle $C_{0}$ in $BH^{0}_{2}-F^{0}+f_{0}+f_{1}$ that contains $e_{0}$. Clearly, $C_{0}$ must contain $f_{0}$, let $H_{0}=C_{0}-f_{0}$, so we denote $C_{0}$ by $\left \langle a^{0},H_{0},b^{0}, a^{0}\right \rangle $. Let $b^{0}a^{3}\in D_{2}-F_{2}$ and $a^{j}b^{j+1}\in D_{2}-F_{2}$ for $j=0,1,2$. Since $\left | F^{j} \right |=0$, by Theorem \ref{2n-2}, there is a Hamiltonian path $H_{j}$ in $BH^{j}_{2}-F^{j}$ that join $a^{j}$ and $b^{j}$ for $j=1,2,3$. Then $\left \langle a^{0},b^{1},H_{1},a^{1},b^{2},H_{2},a^{2},b^{3},H_{3},a^{3},b^{0},H_{0},a^{0}\right \rangle $ (see Fig. \ref{case22}) is a fault-free Hamiltonian cycle in $BH_{3}-F$.

\noindent {\bf Subcase 2.2.} $\left | F^{0} \right |\le 4$.

\vskip 0.02 in

Since $\sum_{j=0}^{3}\left | F^{j} \right |=\left | F \right |-\left | F_{2} \right |  \le 5$ and $\left | F^{0} \right |\le 4$, there exists at most one pivot vertex $a^{0}$ in $BH_{2}^{0}$. Then $3\le \left | F^{0} \right | \le 4$ and $\left | F_{2} \right | \le 5$. Assume that $f_{0}=a^{0}b^{0}$ is an $r$-edge in $e(a^{0})\cap F^{0}$. By Theorem \ref{4n-5}, there is a Hamiltonian cycle $C_{0}$ in $BH^{0}_{2}-F^{0}+f_{0}$. And $C_{0}$ must contain $f_{0}$. Let $H_{0}=C_{0}-f_{0}$ and let $b^{0}a^{3}\in D_{2}-F_{2}$, $a^{j}b^{j+1}\in D_{2}-F_{2}$ for $j=0,1,2$. Since $\left | F^{j} \right |\le 8-3-3=2$, by Theorem \ref{2n-2}, there is a Hamiltonian path $H_{j}$ in $BH^{j}_{2}-F^{j}$ that join $a^{j}$ and $b^{j}$ for $j=1,2,3$. Then $\left \langle a^{0},b^{1},H_{1},a^{1},b^{2},H_{2},a^{2},b^{3},H_{3},a^{3},b^{0},H_{0},a^{0}\right \rangle $ is a fault-free Hamiltonian cycle in $BH_{3}-F$.

\begin{figure}[h]
\centering
\includegraphics[width=65mm]{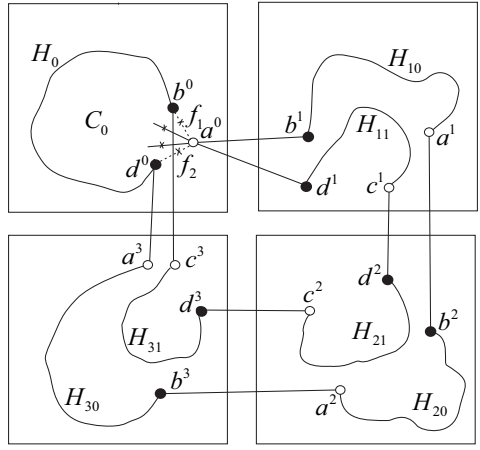}
\caption{Illustration for Case 3.} \label{case3}
\end{figure}

\vskip 0.02 in

\noindent {\bf Case 3.} $\delta (BH^{j}_{2}-F^{j})=0$ for some $j\in \left \{0, 1, 2, 3\right \}$.

\vskip 0.02 in

Since $\sum_{j=0}^{3}\left | F^{j} \right |\le 5$, there is at most one isolated vertex in $BH_{2}^{j}$ for $j=0, 1, 2, 3$. By assumption the isolated vertex $a^{0}\in V(BH^{0}_{2})$. Let $a^{0}b^{1}$, $a^{0}d^{1}\in D_{2}-F_{2}$. Since $\left | F_{2} \right |\le4$, there exist two $r$-edges $f_{1}=a^{0}b^{0}$ and $f_{2}=a^{0}d^{0}$ in $e(a^{0})\cap F^{0}$ such that $c^{3}b^{0}$, $d^{0}a^{3}\in D_{2}-F_{2}$ satisfing $a^{3}\ne c^{3}$. Since $\delta (BH^{0}_{2}-F^{0}+f_{1}+f_{2})=2$ and $\left | F^{0}-f_{1}-f_{2} \right |\le 3$, by Theorem \ref{4n-5}, there is a Hamiltonian cycle $C_{0}$ of $BH^{0}_{2}-F^{0}+f_{1}+f_{2}$ containing $f_{1}$ and $f_{2}$. Let $H_{0}=C_{0}-f_{1}-f_{2}$. Then we can choose $a^{j}b^{j+1}$, $c^{j}d^{j+1}\in D_{2}-F_{2}$ for $j=1,2$. Since $\left | F^{j} \right |\le 8-3-4=1$, by Theorem \ref{2n-3}, there are two vertex-disjoint paths $H_{j0}$ and $H_{j1}$ in $BH^{j}_{2}-F^{j}$ such that $H_{j0}$ joins $a^{j}$ and $b^{j}$ and $H_{j1}$ joins $c^{j}$ and $d^{j}$ for $i=1,2,3$. Hence, $\langle a^{0},b^{1},H_{10},a^{1},b^{2},H_{20},a^{2},b^{3},H_{30},a^{3},d^{0},H_{0},
b^{0},c^{3},H_{31},d^{3},c^{2},H_{21},d^{2},c^{1},H_{11},d^{1},a^{0}\rangle $\\(see Fig. \ref{case3}) is a fault-free Hamiltonian cycle in $BH_{3}-F$.

 }\qed
\end{lemma}

\begin{theorem}{\label{5n-7} \bf .}
Let $F$ be an edge subset of $BH_{n}$ with $\lvert F \rvert\le5n-7$ for $n\ge2$. $BH_{n}-F$ is Hamiltonian if there exists no $f_{4}$-cycles in $BH_{n}-F$ and $\delta (BH_{n}-F)\ge 2$.

\vskip 0.05 in

\noindent {\bf Proof.} {\rm Clearly, it suffices to consider the condition that $\lvert F \rvert=5n-7$. We prove the theorem by induction on $n$. By Theorem \ref{4n-5}, the assertion holds when $n=2$. And it follows from Lemma \ref{8} that the theorem holds when $n=3$. Thus, we assume that the statement holds for $n-1$ wherever $n\ge 4$. Next we consider $BH_{n}$. Since $\sum_{i=0}^{n-1}\lvert F_{i} \rvert = \lvert F \rvert=5n-7$, there is an integer $i\in \left \{0, 1, 2,\dots,n-1 \right \} $ such that $\left | F_{i} \right |\ge 4$. Without loss of generality, we assume that $\left | F_{n-1} \right |\ge 4$. In the following, we decompose $BH_{n}$ into
four subcubes $BH_{n-1}^{j}$ by deleting $D_{n-1}$, where $j=0, 1, 2, 3$. We distinguish the following three cases.
\vskip 0.02 in

\noindent {\bf Case 1.} $\delta (BH^{j}_{n-1}-F^{j})\ge 2$ for all $j\in \left \{0, 1, 2, 3\right \}$.

\vskip 0.02 in

\noindent {\bf Subcase 1.1.} There exists no $f_{4}$-cycles in $BH^{j}_{n-1}-F^{j}$ for all $j\in \left \{0, 1, 2, 3\right \}$.

\vskip 0.02 in

\noindent {\bf Subcase 1.1.1.} $\left | F^{0} \right |=5n-11$.

\vskip 0.02 in

Clearly, $\left | F_{n-1} \right |=4$ and $\left | F^{1} \right |=\left | F^{2} \right |=\left | F^{3} \right |=0$. Since $\left | F_{n-1} \right |=4$ and $\delta (BH^{0}_{n-1}-F^{0})\ge 2$, there exist at least $5n-11-2(2n-4)=n-3$ $r$-edges in $F^{0}$ for $n\ge 4$. Let $f_{0}=a^{0}b^{0}$ be an $r$-edge in $F^{0}$. By induction hypothesis, there is a Hamiltonian cycle $C_{0}$ in $BH^{0}_{n-1}-F^{0}+f_{0}$. If $f_{0}\notin E(C_{0})$, then there must exist an $r$-edge $f_{1}$ on $C_{0}$. We may assume that $f_{0}\in E(C_{0})$. Let $H_{0}=C_{0}-f_{0}$. Then let $a^{0}b^{1}$, $a^{3}b^{0}$, $a^{1}b^{2}$, $a^{2}b^{3}$ be $(n-1)$-dimensional nonfaulty edges. By Theorem \ref{2n-2}, there is a Hamiltonian path $H_{j}$ in $BH^{j}_{2}-F^{j}$ that join $a^{j}$ and $b^{j}$ for $j=1,2,3$. Then $\left \langle a^{0},b^{1},H_{1},a^{1},b^{2},H_{2},a^{2},b^{3},H_{3},a^{3},b^{0},H_{0},a^{0}\right \rangle $ is a fault-free Hamiltonian cycle in $BH_{n}-F$.

\vskip 0.02 in

\noindent {\bf Subcase 1.1.2.} $\left | F^{0} \right | \le 5n-12$.

\vskip 0.02 in

By induction hypothesis, there is a Hamiltonian cycle $C_{0}$ in $BH^{0}_{n-1}-F^{0}$. Since $\left | F_{n-1} \right | \le 5n-7$, there exist at most $\left \lfloor \frac{5n-7}{2} \right \rfloor$ vertices in $V(C_{0})$ incident with two $(n-1)$-dimensional faulty edges. So there exists an $r$-edge $e_{0}=a^{0}b^{0}$ in $E(C_{0})$. Let $a^{0}b^{1}$, $a^{3}b^{0}\in D_{n-1}-F_{n-1}$. Since $\sum_{j=0}^{3}\left | F^{j} \right |=\left | F \right |-\left | F_{n-1} \right |\le5n-11$, $\left | F^{j} \right |\le 4n-9$ for all $j=1,2,3$ and there exists an integer $j\in \left \{ 1,2,3 \right \}$ such that $\left | F^{j} \right | \le 2n-4$.

\vskip 0.02 in

(1) $\left | F^{2} \right | \le 2n-4$. If there exists an $r$-edge $e_{1}=a^{1}b^{1}\in e(b^{1})\cap E(BH^{1}_{n-1})$, let $a^{1}b^{2}\in D_{n-1}-F_{n-1}$. And if there exists an $r$-edge $e_{3}=a^{3}b^{3}\in e(a^{3})\cap E(BH^{3}_{n-1})$, let $a^{2}b^{3}\in D_{n-1}-F_{n-1}$. Since $\left | F^{j} \right | \le 4n-9$ and $\delta (BH^{j}_{n-1}-F^{j})\ge 2$, by Theorem \ref{4n-5}, there is a fault-free Hamiltonian cycle $C_{j}$ in $BH_{n-1}^{j}-F^{j}$ that contains $e_{j}$ for $j=1,3$. Let $H_{j}=C_{j}-e_{j}$ for $j=1,3$. Since $\left | F^{2} \right |  \le 2n-4$, by Theorem \ref{2n-2}, there is a fault-free Hamiltonian path $H_{2}$ in $BH^{2}_{n-1}-F^{2}$ that joins $a^{2}$ and $b^{2}$. Thus, $\left \langle a^{0}, b^{1}, H_{1}, a^{1}, b^{2}, H_{2}, a^{2}, b^{3}, H_{3}, a^{3}, b^{0}, H_{0}, a^{0}\right \rangle $ is a fault-free Hamiltonian cycle in $BH_{n}-F$.

\vskip 0.02 in

(2) $\left | F^{3} \right | \le 2n-4$(the case $\left | F^{1} \right | \le 2n-4$ is similar). If there exists an $r$-edge $e_{1}=a^{1}b^{1}\in e(b^{1})\cap E(BH^{1}_{n-1})$, let $a^{1}b^{2}\in D_{n-1}-F_{n-1}$. And if there exists an $r$-edge $e_{2}=a^{2}b^{2}\in e(b^{2})\cap E(BH^{2}_{n-1})$, let $a^{2}b^{3}\in D_{n-1}-F_{n-1}$. By Theorem \ref{4n-5}, there is a fault-free Hamiltonian cycle $C_{j}$ in $BH_{n-1}^{j}-F^{j}$ that contains $e_{j}$ for $j=1,2$. Let $H_{j}=C_{j}-e_{j}$ for $j=1,2$. Since $\left | F^{3} \right |  \le 2n-4$, by Theorem \ref{2n-2}, there is a fault-free Hamiltonian path $H_{3}$ in $BH^{3}_{n-1}-F^{3}$ that joins $a^{3}$ and $b^{3}$. Thus, $\left \langle a^{0}, b^{1}, H_{1}, a^{1}, b^{2}, H_{2}, a^{2}, b^{3}, H_{3}, a^{3}, b^{0}, H_{0}, a^{0}\right \rangle $ is a fault-free Hamiltonian cycle in $BH_{n}-F$.

\vskip 0.02 in

If one of the aforementioned edges $e_{1}$, $e_{2}$ and $e_{3}$ does not exist, say $e_{1}$, it means that all edges in $e(b^{1})\cap E(BH^{1}_{n-1})$ are non-$r$-edges, then there exist $2(2n-2)=4n-4$ $(n-1)$-dimensional faulty edges adjacent to the edges in $e(b^{1})\cap E(BH^{1}_{n-1})$. It is easy to choose four edges $a^{0}b^{1}$, $a^{1}b^{2}$, $a^{2}b^{3}$, $a^{3}b^{0}\in D_{n-1}-F_{n-1}$. Since $\left | F^{j} \right | \le 5n-7-(4n-4)=n-3\le 2n-4$, by Theorem \ref{2n-2}, there is a fault-free Hamiltonian path $H_{j}$ in $BH^{j}_{n-1}-F^{j}$ that joins $a^{j}$ and $b^{j}$ for  $j=0, 1, 2, 3$. Hence, $\left \langle a^{0}, b^{1}, H_{1}, a^{1}, b^{2}, H_{2}, a^{2}, b^{3},
H_{3}, a^{3}, b^{0}, H_{0}, a^{0}\right \rangle $ is a fault-free Hamiltonian cycle in $BH_{n}-F$.

\begin{figure}
\begin{minipage}[t]{0.5\linewidth}
\centering
\includegraphics[width=65mm]{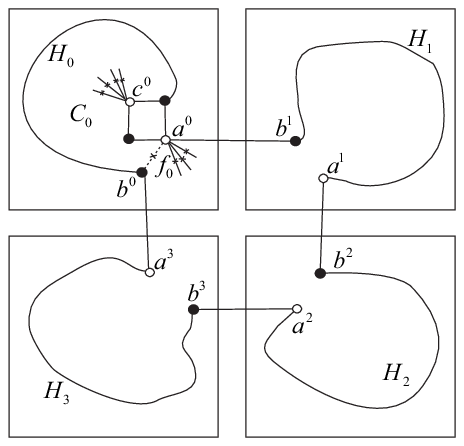}
\caption{Illustration for Subcase 1.2.} \label{ncase2}
\end{minipage}
\begin{minipage}[t]{0.5\linewidth}
\centering
\includegraphics[width=65mm]{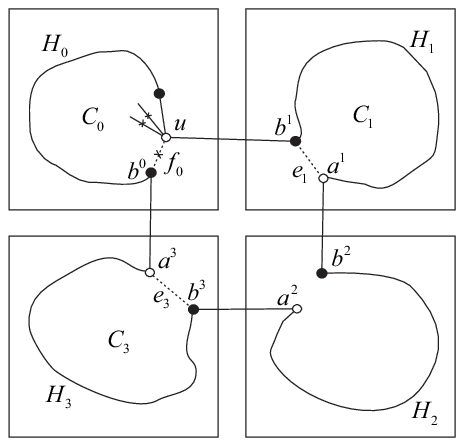}
\caption{Illustration for Subcase 2.1.} \label{ncase31b}
\end{minipage}
\end{figure}

\vskip 0.02 in

\noindent {\bf Subcase 1.2.} There exists an $f_{4}$-cycle in $BH^{j}_{n-1}-F^{j}$ for some $j\in \left \{0, 1, 2, 3\right \}$.

\vskip 0.02 in

An $f_{4}$-cycle is incident with $2(2n-4)=4n-8$ faulty edges, so there exists at most one $f_{4}$-cycle in $BH_{n-1}^{j}-F^{j}$. By our assumption, it follows that the $f_{4}$-cycle lies in $BH_{n-1}^{0}-F^{0}$. Let $a^{0}$, $c^{0}$ be the nonadjacent vertices whose degrees are both two in the $f_{4}$-cycle. Since there exists no $f_{4}$-cycle in $BH_{n}-F$, $a^{0}$ or $c^{0}$, say $a^{0}$, is incident with at least one $(n-1)$-dimensional nonfaulty edge $a^{0}b^{1}\in D_{n-1}-F_{n-1}$. Since $\left | F_{n-1} \right |\le5n-7-(4n-8)=n+1<2(2n-4)$ whenever $n\ge 4$, there exists a faulty $r$-edge $f_{0}=a^{0}b^{0}$ in $BH^{0}_{n-1}$, then let $a^{3}b^{0}\in D_{n-1}-F_{n-1}$. Clearly, there exists no $f_{4}$-cycles in $BH^{0}_{n-1}-F^{0}+f_{0}$, and $\left | F^{0}-f_{0} \right |\le 5n-12$, by induction hypothesis, there is a Hamiltonian cycle $C_{0}$ in $BH^{0}_{n-1}-F^{0}+f_{0}$. Notice that the cycle $C_{0}$ must contain $f_{0}$. Let $H_{0}=C_{0}-f_{0}$ and let $a^{1}b^{2}$, $a^{2}b^{3}\in D_{n-1}-F_{n-1}$. Since $\left | F^{j} \right |\le 5n-7-4-(4n-8)=n-3\le 2n-4$, by Theorem \ref{2n-2}, there is a fault-free Hamiltonian path $H_{j}$ in $BH^{j}_{n-1}-F^{j}$ that joins $a^{j}$ and $b^{j}$ for  $j=1, 2, 3$. Thus, $\left \langle a^{0}, b^{1}, H_{1}, a^{1}, b^{2}, H_{2}, a^{2}, b^{3}, H_{3}, a^{3}, b^{0}, H_{0}, a^{0}\right \rangle $ (see Fig. \ref{ncase2}) is a fault-free Hamiltonian cycle in $BH_{n}-F$.

\vskip 0.02 in

\noindent {\bf Case 2.} $\delta (BH^{j}_{n-1}-F^{j})=1$ for some $j\in \left \{0, 1, 2, 3\right \}$.

\vskip 0.02 in

Since $\left | F \right |-\left | F_{n-1} \right |  \le 5n-7-4=5n-11$ and a pivot vertex in $BH_{n}-D_{n-1}$ is incident with $2n-3$ faulty edges in $F-F_{n-1}$, then there exist at most two pivot vertices in $BH_{n}-D_{n-1}$.

\vskip 0.02 in

\noindent {\bf Subcase 2.1.} There is exactly one pivot vertex in $BH_{n-1}^{j}-F^{j}$ for some $j\in \left \{0, 1, 2, 3\right \}$.

\vskip 0.02 in

Clearly, there exists no $f_{4}$-cycle in each $BH_{n-1}^{j}-F^{j}$ since $2n-3+4+2(2n-4)=6n-7\ge 5n-7$. We may assume that the pivot vertex $u\in V(BH^{0}_{n-1})$. Since $\delta (BH_{n}-F)\ge 2$, $u$ is incident with at least one $(n-1)$-dimensional nonfaulty edge, say $ub^{1}\in D_{n-1}-F_{n-1}$. Since $2n-3+2(2n-3)=6n-9\ge 5n-7$ whenever $n\ge 4$, we can choose an $r$-edge $f_{0}=ub^{0}\in e(u)\cap F^{0}$. If there exists an $f_{4}$-cycle in $BH^{0}_{n-1}-F^{0}+f_{0}$, then $u$ must be a vertex in the $f_{4}$-cycle. Since $\left | F^{0} \right |\le 5n-11$ and $u$ is a pivot vertex, there exists at most one vertex with degree 2 in $BH^{0}_{n-1}-F^{0}$, which is the vertex on the $f_{4}$-cycle. Since $\left | F_{n-1} \right |\le 5n-7-(2n-3)-(2n-4)=n<2(2n-3-1)$, there exists an $r$-edge $f_{1}\in e(u)\cap F^{0}-f_{0}$ and there exists no $f_{4}$-cycles in $BH^{0}_{n-1}-F^{0}+f_{1}$. So, we may assume that there exists no $f_{4}$-cycles in $BH^{0}_{n-1}-F^{0}+f_{0}$. Let $a^{3}b^{0}\in D_{n-1}-F_{n-1}$. By induction hypothesis, there is a Hamiltonian cycle $C_{0}$ in $BH^{0}_{n-1}-F^{0}+f_{0}$. Clearly, the cycle $C_{0}$ must contain $f_{0}$, let $H_{0}=C_{0}-f_{0}$. Notice that $\left | F^{j} \right |\le\sum_{j=1}^{3}\left | F^{j} \right |\le5n-7-(2n-3)-4=3n-8\le4n-9$ for $j=1,2,3$, and there exists an integer $j\in\left \{ 1,2,3 \right \}$ such that $\left | F^{j} \right |\le 2n-4$.

\vskip 0.02 in

(1) $\left | F^{2} \right | \le 2n-4$. Since $\left | F_{n-1} \right |\le 5n-7-(2n-3)=3n-4<2(2n-2)$, there exist an $r$-edge $e_{1}=a^{1}b^{1}\in e(b^{1})\cap E(BH^{1}_{n-1})$ and an $r$-edge $e_{3}=a^{3}b^{3}\in e(a^{3})\cap E(BH^{3}_{n-1})$. Assume that $a^{1}b^{2}$, $a^{2}b^{3}\in D_{n-1}-F_{n-1}$. By Theorem \ref{4n-5}, there is a fault-free Hamiltonian cycle $C_{j}$ in $BH_{n-1}^{j}-F^{j}$ that contains $e_{j}$ for $j=1,3$. Let $H_{j}=C_{j}-e_{j}$ for $j=1,3$. By theorem 4, there is a fault-free Hamiltonian path $H_{2}$ in $BH^{2}_{n-1}-F^{2}$ that joins $a^{2}$ and $b^{2}$. Thus, $\left \langle u, b^{1}, H_{1}, a^{1}, b^{2}, H_{2}, a^{2}, b^{3}, H_{3}, a^{3}, b^{0}, H_{0}, u\right \rangle $ (see Fig. \ref{ncase31b}) is a fault-free Hamiltonian cycle in $BH_{n}-F$.

\vskip 0.02 in

(2) $\left | F^{3} \right | \le 2n-4$(the case $\left | F^{1} \right | \le 2n-4$ is similar). Since $\left | F_{n-1} \right |\le 5n-7-(2n-3)=3n-4<2(2n-2)$, there exists an $r$-edge $e_{1}=a^{1}b^{1}\in e(b^{1})\cap E(BH^{1}_{n-1})$. Assume that $a^{1}b^{2}\in D_{n-1}-F_{n-1}$. And there exists an $r$-edge $e_{2}=a^{2}b^{2}\in e(b^{2})\cap E(BH^{2}_{n-1})$. In addition, assume that $a^{2}b^{3}\in D_{n-1}-F_{n-1}$. By Theorem \ref{4n-5}, there is a fault-free Hamiltonian cycle $C_{j}$ in $BH_{2}^{j}-F^{j}$ that contains $e_{j}$ for $j=1,2$. Let $H_{j}=C_{j}-e_{j}$ for $j=1,2$. Since $\left | F^{3} \right |  \le 2n-4$, by Theorem \ref{2n-2}, there is a fault-free Hamiltonian path $H_{3}$ in $BH^{3}_{n-1}-F^{3}$ that joins $a^{3}$ and $b^{3}$. Thus, $\left \langle u, b^{1}, H_{1}, a^{1}, b^{2}, H_{2}, a^{2}, b^{3}, H_{3}, a^{3}, b^{0}, H_{0}, u\right \rangle $ is a fault-free Hamiltonian cycle in $BH_{n}-F$.

\vskip 0.02 in

\noindent {\bf Subcase 2.2.} There are two pivot vertices $u$ and $v$ in $BH_{n-1}^{j}-F^{j}$ for some $j\in \left \{0, 1, 2, 3\right \}$.

\vskip 0.02 in

If one of the following three conditions (1), (2) and (3) holds, then there exists an integer $i\in \left \{ 0,1,\dots,n-2 \right \}$ such that $u$ and $v$ are both incident with two $i$-dimensional faulty edges and $\left | F_{i} \right |\ge 4$.

\vskip 0.02 in
(1) $u$ and $v$ are not adjacent.
\vskip 0.02 in
(2) $uv$ is a nonfaulty edge.
\vskip 0.02 in
(3) $uv$ is an $(n-1)$-dimensional faulty edge.
\vskip 0.02 in

We may assume that $i=n-2$, then we can decompose $BH_{n}$ into four subcubes $BH_{n-1}^{j,n-2}$ by deleting $D_{n-2}$, where $j=0,1,2,3$. In $BH_{n-1}^{j,n-2}$ for $j=0,1,2,3$, $u$ and $v$ are two vertices whose degrees are 2 or 3, notice that there exists no isolated vertex and at most one pivot vertex in $BH_{n}-D_{n-2}$. Thus, the case is degenerated to the Case 1 or Subcase 2.1.

\vskip 0.02 in

 If $uv$ is a faulty edge and $uv$ is not an $(n-1)$-dimensional edge, we may assume that $u$, $v\in V(BH_{n-1}^{0})$. Notice that $uv$ is an $r$-edge, assume that $ub^{1}$, $a^{3}v$ are $(n-1)$-dimensional nonfaulty edges. By induction hypothesis, there is a Hamiltonian cycle $C_{0}$ in $BH^{0}_{n-1}-F^{0}+f_{0}$. Clearly, the cycle $C_{0}$ must contain $f_{0}$, let $H_{0}=C_{0}-f_{0}$. Let $a^{1}b^{2}$, $a^{2}b^{3}\in D_{n-1}-F_{n-1}$. Since $\left | F^{j} \right |\le 5n-7-4-2(2n-3)+1=n-4<2n-4$, by Theorem \ref{2n-2}, there is a fault-free Hamiltonian path $H_{j}$ in $BH^{j}_{n-1}-F^{j}$ that joins $a^{j}$ and $b^{j}$ for $j=1,2,3$. Therefore, $\left \langle u, b^{1}, H_{1}, a^{1}, b^{2}, H_{2}, a^{2}, b^{3}, H_{3}, a^{3}, v, H_{0}, u\right \rangle $ (see Fig. \ref{ncase32}) is a fault-free Hamiltonian cycle in $BH_{n}-F$.

\begin{figure}
\begin{minipage}[t]{0.5\linewidth}
\centering
\includegraphics[width=65mm]{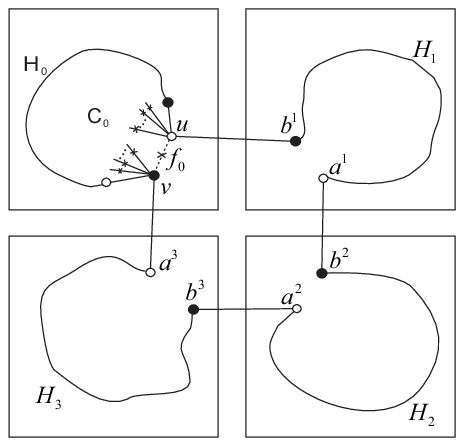}
\caption{Illustration for Subcase 2.2.} \label{ncase32}
\end{minipage}
\begin{minipage}[t]{0.5\linewidth}
\centering
\includegraphics[width=65mm]{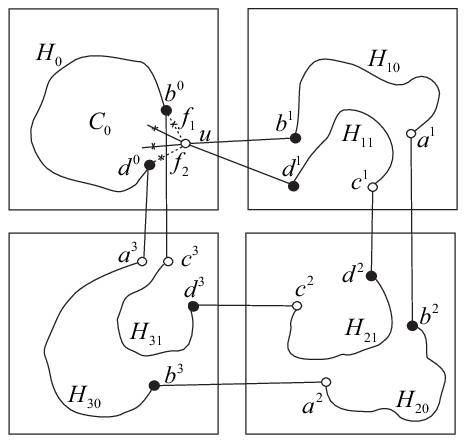}
\caption{Illustration for Subcase 3.1.} \label{ncase41}
\end{minipage}
\end{figure}

\vskip 0.02 in

\noindent {\bf Case 3.} $\delta (BH^{j}_{n-1}-F^{j})=0$ for some $j\in \left \{0, 1, 2, 3\right \}$.

\vskip 0.02 in

Since $\left | F \right |-\left | F_{n-1} \right |  \le 5n-7-4=5n-11$, there exist at most two isolated vertices or at most one pivot vertex and one isolated vertex in $BH_{n-1}^{j}-F^{j}$ for some $j\in \left \{0, 1, 2, 3\right \}$.

\vskip 0.02 in

\noindent {\bf Subcase 3.1.} There is exactly one isolated vertex in $BH_{n-1}^{j}-F^{j}$ for some $j\in \left \{0, 1, 2, 3\right \}$.

\vskip 0.02 in

We may assume that the isolated vertex $u\in BH^{0}_{n-1}$, then $u$ is incident with $2n-2$ edges of $F^{0}$ and two $(n-1)$-dimensional nonfaulty edges $ub^{1}$ and $ud^{1}$.

\vskip 0.02 in

If $\left | F_{i} \right |\le 3$ for each $i\in\left \{ 0,1,\dots,n-2 \right \}$, then $\left | F^{j} \right |\le 3(n-1)-(2n-2)=n-1\le 2n-5$ for $j=1,2,3$. Since $\left | F_{n-1} \right |\le 5n-7-(2n-2)=3n-5<2(2n-2)-3=4n-7$, there exist two $r$-edges $f_{1}=ub^{0}$ and $f_{2}=ud^{0}$ in $e(u)\cap F^{0}$ satisfying $a^{3}b^{0}$, $c^{3}d^{0}\in D_{n-1}-F_{n-1}$ and $a^{3}\ne c^{3}$. Since $\left | F^{0} \right |\le 3(n-1)$ and $u$ is an isolated vertex, there exists no vertices with degree at most two in $BH^{0}_{n-1}-F^{0}$. So, there exists no $f_{4}$-cycle in $BH^{0}_{n-1}-F^{0}+f_{1}+f_{2}$. Since $\left | F^{0}-f_{1}-f_{2} \right |\le 3(n-1)-2=3n-5\le 5n-12$, by induction hypothesis, there is a Hamiltonian cycle $C_{0}$ in $BH^{0}_{n-1}-F^{0}+f_{1}+f_{2}$. Clearly, the cycle $C_{0}$ must contain $f_{1}$ and $f_{2}$, let $H_{0}=C_{0}-f_{1}-f_{2}$. And let $a^{j}b^{j+1}$, $c^{j}d^{j+1}\in D_{n-1}-F_{n-1}$ for $j=1,2$. By Theorem \ref{2n-3}, there are two vertex-disjoint paths $H_{j0}$ and $H_{j1}$ in $BH_{n-1}^{j}-F^{j}$ such that $H_{j0}$ joins $a^{j}$ and $b^{j}$ and $H_{j1}$ joins $c^{j}$ and $d^{j}$ for $j=1,2,3$. Then $\left \langle u,b^{1},H_{10},a^{1},b^{2},H_{20},a^{2},b^{3},H_{30},a^{3},d^{0},H_{0},b^{0},c^{3},H_{31}
,d^{3},c^{2},H_{21},d^{2},c^{1},H_{11},d^{1},u\right \rangle $ \\(see Fig. \ref{ncase41}) is  a fault-free Hamiltonian cycle in $BH_{n}-F$.

\vskip 0.02 in

If there exists an integer $i\in\left \{ 0,1,\dots,n-2 \right \}$ such that $\left | F_{i} \right |\ge 4$, say $\left | F_{n-2} \right |\ge 4$, then we decompose $BH_{n}$ into four subcubes $BH_{n-1}^{j,n-2}$ by deleting $D_{n-2}$, where $j=0,1,2,3$. Notice that $u$ is a vertex whose degree is 2 in $BH_{n-1}^{j,n-2}$ for some $j\in \left \{0, 1, 2, 3\right \}$. If there exists no isolated  vertex in $BH_{n-1}^{j,n-2}$ for each $j\in \left \{0, 1, 2, 3\right \}$, then this case  is degenerated to Case 1 or Case 2. Otherwise, if there still exists an isolated vertex $v$ in $BH_{n-1}^{j,n-2}$ for some $j\in \left \{0, 1, 2, 3\right \}$, then there exists an integer $i\in \left \{ 0,1,\dots,n-3 \right \}$ such that $u$ and $v$ are both incident with two $i$-dimensional faulty edges and $\left | F_{i} \right |\ge 4$, we may assume that $i=n-3$. Then we decompose $BH_{n}$ into four subcubes $BH_{n-1}^{j,n-3}$ by deleting $D_{n-3}$, where $j=0,1,2,3$. Notice that there exists no isolated  vertex in $BH_{n-1}^{j,n-3}$ for $j\in \left \{0, 1, 2, 3\right \}$ as $5n-7-(2n-2)-(2n-2)=n-3$, the case is degenerated to Case 1 or Case 2.

\vskip 0.02 in

\noindent {\bf Subcase 3.2.} There are two isolated vertices $u$ and $v$ or a pivot vertex $u$ and an isolated vertex $v$ in $BH_{n-1}^{j}-F^{j}$ for some $j\in \left \{0, 1, 2, 3\right \}$.

\vskip 0.02 in

It can be deduced that there exists an integer $i\in \left \{ 0,1,\dots,n-2 \right \}$ such that $u$ and $v$ are both incident with two $i$-dimensional faulty edges and $\left | F_{i} \right |\ge 4$, we may assume that $i=n-2$. Then we decompose $BH_{n}$ into four subcubes $BH_{n-1}^{j,n-2}$ by deleting $D_{n-2}$, where $j=0,1,2,3$. Since $5n-7-(2n-2)-(2n-3)\ge n-2$, there exists no isolated  vertex in $BH_{n-1}^{j,n-2}$ for $j\in \left \{0, 1, 2, 3\right \}$, the case is degenerated to Case 1 or Case 2.

}\qed
\end{theorem}

\section{Conclusion}
In this paper, we show that there exists a fault-free Hamiltonian cycle in $BH_{n}-F$ for $n\geq 2$ with $\left | F \right |\leq 5n-7$ if the degree of every vertex in $BH_{n}-F$ is at least two and there exists no $f_{4}$-cycles in $BH_{n}-F$, which improves the upper bound of the number of faulty edges can be tolerated in $BH_n$. It is interesting to further consider Hamiltonicity of the faulty $BH_{n}$ that contains no $f_{6}$-cycle.

\end{document}